\documentclass{amsart}
\usepackage{amssymb}
\usepackage[all]{xy}
\newcommand\datver[1]{\def\datverp%
 {\par\boxed{\boxed{\text{#1; Run: \today}}}}}
\datver{1.0J; Revised: 19-11-2003}

\newcommand\boxb[1]{\square_b}

\newcommand\ff{\operatorname{ff}}

\numberwithin{equation}{section}
\renewcommand{\theequation}{\arabic{equation}}
\let\Appendix\appendix
\newcommand\paperbody%
        {\renewcommand{\theequation}{\thesection.\arabic{equation}}}
\renewcommand\appendix{\Appendix%
\renewcommand{\theequation}{A.\arabic{equation}}}
\newtheorem{lemma}{Lemma}
\newtheorem{proposition}{Proposition}
\newtheorem{corollary}{Corollary}
\newtheorem{theorem}{Theorem}
\newtheorem{non-theorem}{Non-Theorem}

\theoremstyle{remark}

\newtheorem{remark}{Remark}

\newcommand\Kc{K_c}

\newcommand\bo{\operatorname{b}}

\newcommand\cusp{\operatorname{cu}}
\newcommand\scusp{\operatorname{s-cu}}
\newcommand\nd{\operatorname{nd}}

\newcommand\cFTs{{}^{\Phi}\overline{T}\kern-1pt{}^*}

\newcommand\sus{\operatorname{sus}}

\newcommand\even{\text{even}}
\newcommand\odd{\text{odd}}

\newcommand\Tr{\operatorname{Tr}}

\newcommand\bTr{\overline{\operatorname{Tr}}}

\newcommand\RTr[1]{\operatorname{Tr_{R,#1}}}

\newcommand\GL{\operatorname{GL}}

\newcommand\tx{\tilde{x}}

\newcommand\com[1]{\overline{#1}}

\newcommand\ie{i\@.e\@. }

\newcommand\bbC{\mathbb C}

\newcommand\bbN{\mathbb N}

\newcommand\bbR{\mathbb R}
\newcommand\bbS{\mathbb S}

\newcommand\bbZ{\mathbb Z}

\newcommand\CIc{{\mathcal{C}}^{\infty}_c}

\newcommand\CI{{\mathcal{C}}^{\infty}}

\newcommand\CmI{{\mathcal{C}}^{-\infty}}

\newcommand\Diag{\operatorname{Diag}}

\newcommand\cFNs{{}^{\Phi}\overline N\kern-1pt{}^*}

\newcommand\ind{\operatorname{ind}}

\newcommand\Hom{\operatorname{Hom}}

\newcommand\Id{\operatorname{Id}}

\newcommand\dCI{\dot{\mathcal{C}}^{\infty}}

\newcommand\nul{\operatorname{null}}

\newcommand\pa{\partial}

\newcommand\supp{\operatorname{supp}}

\newcommand\inn{\operatorname{int}}

\newcommand\cl{\operatorname{cl}}

\renewcommand\Re{\operatorname{Re}}
\renewcommand\Im{\operatorname{Im}}

\newcommand\Mand{\text{ and }}

\newcommand\Mat{\text{ at }}

\newcommand\Min{\text{ in }}

\begin{document}
\title[Families index with boundary]
{Families index for pseudodifferential operators on manifolds with boundary}

\author{Richard Melrose}
\email{rbm@math.mit.edu}
\author{Fr\'ed\'eric Rochon}
\email{rochon@math.mit.edu}
\address{Department of Mathematics, Massachusetts Institute of Technology}
\thanks{RBM acknowledges support from the National Science Foundation
under Grant DMS-0104116, FR from the Natural Sciences and Engineering
Research Council of Canada}
\begin{abstract} An analytic index is defined for a family of cusp
  pseudodifferential operators, $P_b,$ on a fibration with fibres which are
  compact manifolds with boundaries, provided the family is elliptic and has
  invertible indicial family at the boundary. In fact there is always a
  perturbation $Q_b$ by a family of cusp operators of order $-\infty$ such
  that each $P_b+Q_b$ is invertible. Thus any elliptic family of symbols has a 
  realization as an invertible family of cusp pseudodifferential operators,
  which is a form of the cobordism invariance of the index. A crucial role
  is played by the weak contractibility of the group of cusp smoothing
  operators on a compact manifold with non-trivial boundary and the
  associated exact sequence of classifying spaces of odd and even
  K-theory.
\end{abstract}
\maketitle

\tableofcontents

\section*{Introduction}

For a smooth fibration of compact manifolds
\begin{equation}
\xymatrix{Z\ar@{-}\ar[r]&M\ar[d]^{\phi}\\&B}
\label{fipomb.1}\end{equation}
with model fibre $Z,$ with boundary, we show that any elliptic symbol on
the fibres, for bundles $E,$ $F$ over the total space, can be quantised to
an invertible family of cusp pseudodifferential operator. This is the
cobordism invariance of the index in this context, amounting to the fact
that for any elliptic family $P\in\Psi_{\cusp}^m(M/B;E,F)$ the associated
indicial family, which is an elliptic loop $I(P,s)\in\CI(\bbR;\Psi^m(\pa
M/B;E,F)),$ defines the trivial class in $K^1(B).$ The homotopy classes of
Fredholm perturbations $P+Q,$ with $Q\in\Psi_{\cusp}^{-\infty}(M/B;E,F),$
form a model for $K^0(B),$ with the group structure given by a relative index
formula in terms of the indicial families
\begin{equation}
\ind(P+Q_2)-\ind(P+Q_1)=[I(P+Q_2)I(P+Q_1)^{-1}]\in K^0(B).
\label{CnvRcn.11}\end{equation}
The result that proves crucial is the weak contractibility of the group of
the invertible operators of the form $\Id+A$ with $A$ a cusp operator of
order $-\infty.$

To put these results in perspective, recall the index theorem for families
of pseudodifferential operators, of Atiyah and Singer; this may be viewed
as a topological description of the analytic index. For a fibration of
compact manifolds, as in \eqref{fipomb.1} but without boundary, an elliptic
family of (classical) pseudodifferential operators on the fibres, acting
between sections of vector bundles, $P\in\Psi^m(M/B;E,F)$ has an associated
symbol, $\sigma(P)\in\CI(S^*(M/B);E,F),$ which is invertible precisely by
the definition of ellipticity. As a family of operators from $\CI(M;E)$ to
$\CI(M;F),$ as infinite rank bundles over $B,$ $P$ is Fredholm. It may be
perturbed, by an element of $\Psi^{-\infty}(M/B;E,F),$ hence without
changing the principal symbol, to have null space of constant rank, and
this defines a K-class on the base
\begin{equation}
\ind_{\text{a}}(P)=[\nul(P)\ominus\nul(P^*)]\in K^0(B)
\label{CnvRcn.13}\end{equation}
which is independent of the choice of perturbation. In particular in this
boundaryless case, the space of Fredholm perturbations $P+Q$ with
$Q\in\Psi^{-\infty}(M/B;E,F)$ is contractible and indeed the index class
only depends on the symbol $\sigma(P).$ The symbol itself defines an
element of $\Kc^0(T^*(M/B)),$ which may be identified with the relative group
for the radial compactification of the fibres of the fibre cotangent
bundle, in which case the fibre cosphere bundle may be identified with its
boundary:
\begin{equation}
[\sigma(P)]\in\Kc^0(T^*(M/B))=K^0(\com{T}^*(M/B);S^*(M/B)).
\label{CnvRcn.14}\end{equation}
The topological index is defined via a Gysin map 
\begin{equation}
\ind_{\text{t}}:\Kc^0(T^*(M/B))\longrightarrow K^0(B)\Mand
\ind_{\text{t}}[\sigma(P)]=\ind_{\text{a}}(P)
\label{CnvRcn.15}\end{equation}
where the equality of these two notions of index is a principal result of
\cite{Atiyah-Singer1}.

Dirac operators, associated to Hermitian Clifford modules, are an
important special case of the Atiyah-Singer theorem. Bismut and
Cheeger, \cite{MR91e:58180}, \cite{MR91e:58181}, generalised
the families index theorem in this case to manifolds with boundary,
following the theorem of Atiyah, Patodi and Singer in the case of a single
operator. Bismut and Cheeger proceeded under the assumption that the family
of Dirac operators induced on the boundary had null space of constant
rank. Topologically, this is not an unreasonable assumption since it
is stably equivalent to the vanishing of the class in $K^1(B)$ defined
by the (self-adjoint) boundary family; this vanishing is an expression of
the cobordism invariance of the index. However, analytically it is not
reasonable, since these null spaces have no inherent stability, and it was
removed in \cite{MR99a:58144}, 
\cite{MR99a:58145}. In terms of the original context of Atiyah,
Patodi and Singer, and also of Bismut and Cheeger, it was shown there that
for any family of Dirac operators, $\eth,$ associated to a Hermitian
Clifford module $E$ on the fibres of a fibration with
model fibre $Z,$ a (connected) compact manifold with non-trivial
boundary, there is a smooth family of self-adjoint pseudodifferential
projections $\Pi\in \Psi^0(\pa M/B;E)$ such that for each $b\in B,$
$\Pi_b-\Pi_+$ is spectrally finite for the family of boundary Dirac operators
$\eth_0.$ Here $\Pi_+$ is the Atiyah-Patodi-Singer projection onto
the span of the eigensections of $\eth_0$ with positive
eigenvalue. It follows that the family of operators
\begin{equation}
\eth^{+}:\{u\in\CI(M/B;E^+);\Pi(u\big|_{\pa M})=0\}\longrightarrow \CI(M/B;E^-)
\label{CnvRcn.16}\end{equation}
is Fredholm. Thus the Fredholm operator depends on a choice of
projection $\Pi.$ The existence of such a smooth family of projections
is shown in \cite{MR99a:58144} to be equivalent to the (known)
vanishing of the odd index, $[\eth_0]\in K^1(B).$

The analysis of the index in \cite{MR99a:58144} is carried out in
the context of the b-calculus of pseudodifferential operators. This
corresponds to the addition of a `cylindrical end' to the
compact manifold with boundary, the method used by Atiyah, Patodi and
Singer to obtain their formula in the case of a single Dirac
operator. The b-calculus actually corresponds to the asymptotic
$\bbR_+$-invariance of any compact manifold with boundary near its
boundary. That is, rather than considering a `cylindrical end' as a
non-compact manifold it is considered as a compact manifold with
boundary, obtained in effect by the exponential compactification of the
end, $\tx=e^{-t}$ where $t$ is the translation variable on the
end. The effect of this compactification is to replace the fully
elliptic Dirac operator (of product type) by a (family of) degenerate
`b-elliptic' Dirac operator associated to a b-metric, meaning one of the form 
\begin{equation}
g=\frac{d\tx^2}{\tx^2}+h(\tx)
\label{CnvRcn.18}\end{equation}
near the boundary, with $h$ smooth and inducing a metric on the
boundary. Near the boundary then, an associated Dirac operator is elliptic
`as a function of' the tangent vector fields $\tx\pa_{\tx},$ $\pa_{y_j}$
where $y_i$ are boundary coordinates. In the single operator case this is
discussed extensively in \cite{MR96g:58180}.

In the families case in \cite{MR99a:58144,MR99a:58145}, the Fredholm
family for which the index is actually computed is obtained by
perturbation, namely $\eth+Q$ where $Q$ is a b-pseudodifferential operator
of order $-\infty.$ The non-compactness of such a perturbation is captured
by the non-triviality of its indicial family and the perturbation is chosen
precisely so that the indicial family of $\eth+Q$ is invertible. We make a
similar analysis here but in the pseudodifferential case.

There is another compactification of the cylindrical end, which we use here
and which achieves essentially the same effect. Namely one can compactify
by radial inversion, introducing $x=1/t$ where $t$ is the translation
variable along the cylindrical end. Thus $x$ and $\tx$ are related
transcendentally by $x=1/\log(1/\tx).$ The asymptotically
translation-invariant metric on the cylindrical end then takes the `cusp' form
\begin{equation}
g=\frac{dx^2}{x^4}+h(x)
\label{CnvRcn.21}\end{equation}
in place of \eqref{CnvRcn.18}. Just as the (intrinsic) b-structure
on a compact manifold with boundary is associated to an algebra of
pseudodifferential operators, $\Psi_{\bo}^*(Z),$ there is an algebra of
cusp operators $\Psi_{\cusp}^*(Z)$ which however is not quite intrinsic but
depends on the choice of a boundary defining function $x.$ Another choice
$x'$ yields the same structure if $x'=cx+O(x^2)$ where $c>0$ is
constant. This can also be identified as a choice of trivialization of the
normal bundle to the boundary
\begin{equation}
N(\pa Z)\cong \pa Z\times L
\label{CnvRcn.19}\end{equation}
where $L$ is a 1-dimensional real oriented vector space.

Our approach depends crucially on the structure of the cusp algebra; this
is a special case of the class of algebras discussed in
\cite{Mazzeo-Melrose4}. For any bundles $E$ and $F$ over a compact manifold
with boundary it associates a space of operators 
\begin{equation}
\Psi_{\cusp}^m(Z;E,F)\ni A:\CI(Z;E)\longrightarrow \CI(Z;F),
\label{fipomb.2}\end{equation}
where the choice of a cusp
structure is suppressed in the notation. These operators also map
$\dCI(Z;E)$ into $\dCI(Z;F)$ where $\dCI(Z;E)\subset\CI(Z;E)$ is the
subspace of elements vanishing at the boundary in the sense of Taylor
series. The algebra is $*$-invariant, so the elements also act naturally
on the dual spaces, consisting of extendible and supported distributional
sections and restrict to act on the appropriate Sobolev spaces. In
particular
\begin{equation}
P\in\Psi^m_{\cusp}(X;E,F)\Longrightarrow
P:H^M_{\cusp}(X;E)\longrightarrow H^{M-m}_{\cusp}(X;F),\ \forall\ M\in\bbR. 
\label{CnvRcn.20}\end{equation}
Over the interior a cusp pseudodifferential operators reduces to a
pseudodifferential operator in the usual sense. The symbol map is
therefore defined over the interior and it extends by continuity to
define the `cusp' symbol map which gives a well-defined short exact sequence
\begin{equation}
\Psi^{m-1}_{\cusp}(Z;E,F)\longrightarrow
\Psi^m_{\cusp}(Z;E,F)\overset{\sigma_m}
\longrightarrow \CI({}^{\cusp}S^*Z;\hom(E,F)\otimes R^m).
\label{CnvRcn.22}\end{equation}
Here ${}^{\cusp}S^*Z={}^{\cusp}T^*Z\setminus 0/\bbR_+$ is the cusp cosphere
bundle, the intrinsic sphere bundle of the cusp cotangent bundle. The
latter is a bundle canonically associated to the cusp structure. It is
isomorphic to the usual cotangent bundle, but not naturally so. The cusp
tangent bundle, of which ${}^{\cusp}T^*Z$ is the dual, is defined precisely
so that a metric \eqref{CnvRcn.21} is a non-degenerate smooth fibre metric
on it. The trivial line bundle $R^m$ in \eqref{CnvRcn.22} is the bundle
over ${}^{\cusp}S^*Z$ with sections which are functions over
${}^{\cusp}S^*Z\setminus0$ which are positively homogeneous of degree $m.$
The symbol map is multiplicative
\begin{multline}
A\in \Psi^m_{\cusp}(Z;E,F),\ B\in
\Psi^{m'}_{\cusp}(Z;F,G)\Longrightarrow \\
B\circ A\in\Psi^{m+m'}_{\cusp}(Z;E,G),\
\sigma_{m+m}(B\circ A)=\sigma_{m'}(B)\circ\sigma_m(A).
\label{CnvRcn.23}\end{multline}

For (symbolically) elliptic operators the Fredholm condition is captured by the
normal operator, or equivalently the indicial family. For simplicity
consider this initially under the assumption that the boundary is
connected. For a cusp structure on a compact manifold with boundary, the normal
operator is a pseudodifferential operator on the `model space' at the
boundary $L\times\pa Z,$ where $L$ is the oriented one-dimensional real
vector space from \eqref{CnvRcn.19}. Corresponding to the (intrinsic) `asymptotic
translation invariance' of elements of the cusp calculus, it is
invariant under translations in $L;$ in fact its kernel, as a partial
convolution operator, is also rapidly decreasing at infinity in $L.$
The algebra of such pseudodifferential operators is the
`$L$-suspended' algebra, discussed for example in \cite{MR96h:58169} and
described in the appendix. In the setting of a fibration by compact
manifolds with boundary the normal operator corresponds to a short exact
sequence with the corresponding bundle of algebras as image
\begin{equation}
x\Psi^{m}_{\cusp}(M/B;E,F)\longrightarrow\Psi^m_{\cusp}(M/B;E,F)\\
\overset{N}\longrightarrow \Psi_{\sus(L)}^m(\pa M/B;E,F).
\label{CnvRcn.25}\end{equation}
The normal map, $N,$ is multiplicative and if $P$ is an elliptic family the
existence of an inverse $N(P)^{-1}\in \Psi_{\sus(L)}^{-m}(\pa M/B;F,E)$ to
$N(P)$ implies that it is a Fredholm family as a map from the bundle
$\dCI(M;E)$ to $\dCI(M;F)$ over $B.$ For the action on appropriate
Sobolev spaces this condition is necessary and sufficient. In general there
is a normal operator associated to each boundary hypersurface and (see the
Appendix) smoothing terms between them.

The translation-invariance of $N(P)$ allows it to be realized, by
Fourier transform in $L,$ as a 1-parameter family of
pseudodifferential operators on the fibres of $\pa M;$ this is the
indicial family $I(P,s),$ $s\in\bbR.$ Given the ellipticity of $P,$
the invertibility of $N(P)$ in the suspended algebra is equivalent to
the existence of an inverse $I(P,s)^{-1}\in\Psi^{-m}(\pa M/B;F,E)$ for
each $s\in\bbR.$ Note that $I(P,s)$ depends on the trivialization of
$L$ and, once a corresponding boundary defining function has been chosen,
can be defined directly from $P$
\begin{multline}
I(P,s)v=\left(e^{-is/x}P(e^{is/x}\tilde v)\right)\big|_{\pa
  M}\in\CI(\pa M;F),\\ \forall\
  v\in\CI(\pa M;E)\text{ with }\tilde v\in\CI(M;E),\ \tilde v\big|_{\pa M}=v.
\label{CnvRcn.26}\end{multline}
We denote the set of boundary hypersurfaces of a manifold with boundary $X$
by $M_1(X)$ and by $I_{HG}(P,s)$ the indicial family acting between the
boundary hypersurfaces $H$ and $G.$ The full indicial family, $I(P,s)$ is
an $N\times N$ matrix of operators between the $N$ components of the
boundary.

A general elliptic family $P'\in\Psi_{\sus(L)}^m(M'/B';E,F),$ for a
fibration $\phi:M'\longrightarrow B'$ of compact manifolds without
boundary, defines a class in $K^1(B')$ which is precisely the obstruction
to the existence of a perturbation
$Q'\in\Psi_{\sus(L)}^{-\infty}(M'/B';E,F)$ to an invertible family with
$(P'+Q')^{-1}\in\Psi_{\sus(L)}^{-m} (M'/B';F,E).$ The cobordism invariance
of the index in this case takes the form of the vanishing of the $K^1$
class for the boundary family of any elliptic family of cusp
operators. This is shown in Theorem~\ref{5.3.2003.1} below. In particular
there are always Fredholm perturbations of this type and in
Theorem~\ref{5.3.2003.3} we compute the difference of the families index
for two such perturbations and show that the homotopy classes of such
perturbations form a model for $K^0(B).$ As shown in
Theorem~\ref{fipomb.21}, there is an invertible perturbation which defines
the origin in this space.  The analogous `odd' families theorems (in this
case for the corresponding suspended algebras) are also discussed. In a
future publication we shall give a formula for the Chern character of the
families index, involving eta forms on the boundary and analyse the
determinant bundle.

In the first section the index of operators of the form $\Id+A$ where $A$
is a cusp pseudodifferential operator of order $-\infty$ is discussed. The
group, $G^{-\infty}_{\cusp}(X;E),$ of invertible operators of this form is
introduced in Section~2 and is related to the normal subgroup of true
smoothing operators; the corresponding short exact sequence of groups is
used to prove that $G^{-\infty}_{\cusp}(X;E)$ is weakly contractible. This
contractibility is in turn the key to the proof, in Section~3, that any
family of cusp operators which is elliptic has a perturbation by a family
of cusp operators of order $-\infty$ to a Fredholm invertible family; this
is refined in Section~4 to show that there is such a perturbation to an
invertible family. The relative index between Fredholm perturbations is
shown to be an isomorphism to even K-theory in Section~5 and the odd case
is discussed, through suspension, in Section~6. A brief discussion of cusp
pseudodifferential operators can be found in the Appendix.

The authors thank Sergiu Moroianu for a comment on the manuscript.

\paperbody

\section{Reduced index formula}

If $X$ is a compact manifold with boundary and $E$ is a smooth vector
bundle over it, the elements of $\Psi^{-\infty}_{\cusp}(X;E)$ which are
compact form the ideal 
\begin{equation}
x\Psi^{-\infty}_{\cusp}(X;E)\hookrightarrow \Psi^{-\infty}_{\cusp}(X;E)
\overset{N}\longrightarrow \Psi^{-\infty}_{\sus}(\pa X;E)
\label{fipomb.39}\end{equation}
where the quotient is given by the normal operator. The `reduced' formula
computes the index of a Fredholm operator of the form $\Id+A,$
$A\in\Psi^{-\infty}_{\cusp}(X;E).$ It reappears later as the relative index
formula for the full algebra.

\begin{proposition}\label{7.3.2003.1} If
  $A\in\Psi_{\cusp}^{-\infty}(X;E)$ then $\Id+A$ is Fredholm acting on
  $L^2_{\cusp}(X;E)$ (or on any of the spaces $x^rH^m_{\cusp}(X;E),$
  $r,m\in\bbR)$ if and only if the indicial family $\Id+I(A,s)$
  is invertible, acting on $L^2(\pa X;E),$ for all $s\in\bbR,$ and then
  the index can be expressed in terms of the winding number of the
  Fredholm determinant of the indicial family
\begin{equation}
\ind(\Id+A)=\operatorname{wn}(\det(\Id+I(A,s))). 
\label{7.3.2003.2}\end{equation}
\end{proposition}

\begin{proof} The necessity of the condition on the indicial family is part
of the structure theory of the cusp calculus. To prove it, suppose that
$\Id+I(A,s_0)$ is not invertible for some $s_0\in\bbR.$ By the Fredholm
alternative, this implies that $(\Id+I(A,s_0))\phi =0$ for some non-trivial
$\phi \in\CI(\pa X;E)=\bigoplus_{H\in M_1(X)}\CI(H;E).$ This can be used to
construct a family $u_t\in L^2(X;E)$ with $\|u_t\|_{L^2}=1$ which is weakly
convergent to $0$ but is such that $\|(\Id+A)u_t\|\to 0.$ It
follows from this that $\Id+A$ cannot be Fredholm on $L^2.$

If $A\in\Psi_{\cusp}^{-\infty}(X;E)$ is such that $\Id+I(A,s)$ is
invertible for all $s\in\bbR$ then a parametrix for it of the form $\Id+B,$
$B\in\Psi_{\cusp}^{-\infty}(X;E)$ can be constructed. Here
$\Id+I(B,s)=(\Id+I(A,s))^{-1}$ and standard inductive and asymptotic
summation arguments are used. It can be arranged that $\Id+B$ is the
generalised inverse, so has null space exactly the null space of $\Id+A^*$
and range the orthocomplement to $\nul(\Id+A).$

Then Calder\'on's formula for the index becomes 
\begin{equation}
\ind(\Id+A)=\Tr([\Id+A,\Id+B])=\Tr([A,B])=\bTr([A,B]).
\label{7.3.2003.5}\end{equation}
Assuming for simplicity of notation that $\pa X$ is connected we can apply
the trace-defect formula, \eqref{fipomb.33}, and the standard formula for
the logarithmic derivative of the determinant
\begin{multline}
\ind(\Id+A)=\frac1{2\pi i}\int_{\bbR}\Tr\big(I(B,s)\frac{d}{ds}I(A,s)\big)ds\\
=\frac1{2\pi i}\int_{\bbR}\Tr\big((\Id+I(B,s))\frac{d}{ds}I(A,s)\big)ds\\
=\frac1{2\pi i}\int_{\bbR}\Tr\big((\Id+I(A,s))^{-1}\frac{d}{ds}I(A,s)\big)ds\\
=\operatorname{wn}(\det(\Id+I(A,s)).
\label{7.3.2003.6}\end{multline}
\end{proof}

Let $\dot\Psi^{-\infty}(X;E)$ be the space of smoothing operators with
kernels vanishing to infinite order at the boundary.

\begin{proposition}\label{7.3.2003.7} The index map
\begin{equation*}
\ind:\{A\in\Psi_{\cusp}^{-\infty}(X;E);\Id+A\text{ is Fredholm on
  }L^2_{\cusp}(X;E)\} \longrightarrow \bbZ
\end{equation*}
is surjective and the null space is 
\begin{multline}
\{A\in\Psi^{-\infty}_{\cusp}(X;E);\ \exists\
S\in\dot\Psi^{-\infty}(X;E)\text{ with }\\
(\Id+A+S)^{-1}=\Id+B,\ B\in\Psi^{-\infty}_{\cusp}(X;E)\}.
\label{fipomb.38}\end{multline}
\end{proposition}

\begin{proof} The surjectivity of the index follows immediately from
Proposition~\ref{7.3.2003.1} above. If $\Id+A$ has
index zero then its null space is a finite dimensional subspace of
$\dCI(X;E)$ with the same dimension as another subspace of $\dCI(X;E)$
which is a complement to its range. Taking a finite rank operator
$S\in\dot\Psi^{-\infty}(X;E)$ which is an isomorphism between them and is
trivial on the orthocomplement of the null space shows that the space in
\eqref{fipomb.38} contains the set of operators of index zero; equality
follows from the fact that elements of $\dot\Psi^{-\infty}(X;E)$ are
compact on $L^2(X;E)$ so such perturbations do not change the index.
\end{proof}

\section{The cusp-smoothing group}

In what follows, it will be convenient to consider differential forms on
infinite dimensional spaces. This presents no serious problem, with the
forms well-defined multilinear forms on the tangent space at each point but
in any case the reader may, if he prefers, always consider the restriction
of these forms to finite dimensional submanifolds.

The group  
\begin{multline}
G^{-\infty}_{\cusp}(X;E)=\{A\in\Psi^{-\infty}_{\cusp}(X;E);\Id+A\text{ is
  invertible with }\\
(\Id+A)^{-1}=\Id+B,\ B\in\Psi^{-\infty}_{\cusp}(X;E)\}
\label{fipomb.40}\end{multline}
plays a crucial role in the analysis below. In particular, we show in this
section that it is weakly contractible, \ie that all its homotopy groups
are trivial.

First consider the normal subgroup consisting of invertible smoothing
perturbations of the identity.

\begin{proposition}\label{fipomb.6} 
The group of smoothing operators on any compact manifold with boundary 
\begin{equation}
\dot G^{-\infty}(X;E)=\left\{\Id+A,\
A\in\dot\Psi^{-\infty}(X;E);(\Id+A)^{-1}=\Id+B,\
B\in\dot\Psi^{-\infty}(X)\right\}
\label{fipomb.3}\end{equation}
is a classifying group for odd K-theory so for all $k\in \bbN,$ $\pi_{2k}(\dot
G^{-\infty}(X;E))=\{0\},$ and isomorphisms
\begin{equation}
\pi_{2k-1}(\dot G^{-\infty}(X;E))\longrightarrow \bbZ
\label{fipomb.4}\end{equation}
are defined by pull-back and integration of the basic (or index) classes, 
\begin{equation}
\beta_k^{\odd}=c_k^{\odd}\Tr\left((a^{-1}d a)^{2k-1}\right),\
c_k^{\odd}=\frac{1}{(2\pi i)^{k}}\frac{(k-1)!}{(2k-1)!}.
\label{fipomb.5}\end{equation}
\end{proposition}

\begin{proof} Finite rank approximation reduces this to the
computation of the stable homotopy groups of $\GL(N,\bbC)$ and of the index
classes. The explicit constant $c_k^{\odd}$ can be found in Theorem
19.3.1 of the book of H\"ormander \cite{Hormander3} and using section 2 of the
paper of Atiyah, \cite{MR35:2299}.
\end{proof}

\begin{proposition}\label{fipomb.7} For any compact manifold without
  boundary, $Y,$ the group of suspended smoothing operators
\begin{multline}
G^{-\infty}_{\sus}(Y;E)=\big\{\Id+A;A\in\mathcal{S}(\bbR;\Psi^{-\infty}(Y;E)),\\
(\Id+A)^{-1}=\Id+B,\ B\in\mathcal{S}(\bbR;\Psi^{-\infty}(Y;E))\big\}
\label{fipomb.8}\end{multline}
is a classifying group for even K-theory, so
$\pi_{2k-1}(G^{-\infty}_{\sus}(Y;E)=\{0\}$ with the isomorphisms 
\begin{equation}
\pi_{2k}(G^{-\infty}(Y;E))\longrightarrow \bbZ
\label{fipomb.9}\end{equation}
defined by pull-back and integration of the basic classes
\begin{equation}
\beta_k^{\even}=c_k^{\even}\int_{\bbR}\Tr\left((a^{-1}d
a)^{2k}a^{-1}\frac{da}{ds}\right),\
c_k^{\even}=\frac{1}{(2\pi i)^{k+1}}\frac{k!}{(2k)!}=\frac{c_k^{\odd}}{4\pi i}.
\label{fipomb.10}\end{equation}
\end{proposition}

\begin{remark}\label{fipomb.68} It is conventional to suppose that a
  manifold is connected. However this conflicts with the convenient
  assertion that the boundary of a manifold with boundary is itself a
  manifold, rather than the more awkward `finite union of disjoint
  manifolds'. In Propositions~\ref{fipomb.6} and \ref{fipomb.7} we can
  allow the wider, not necessarily connected notion, of a manifold provided
  operators are allowed \emph{between} the components. Thus, for example,
  the space of smoothing operators in this case is isomorphic to 
\begin{equation}
\CI(Y^2)=\bigoplus_{i,j=1}^N\CI(Y_i\times Y_j)
\label{fipomb.69}\end{equation}
where $Y=\bigcup_{i=1}^NY_i$ is the decomposition into components. This is
our convention in this paper.
\end{remark}

\begin{proof} If $f:\bbS^{2k+1}\rightarrow G^{-\infty}(Y;E)$ is a
generator of $\pi_{2k+1}( G^{-\infty}(Y;E))\cong \bbZ,$ it gives rise
in a natural way to a generator of $\pi_{2k}(G^{-\infty}_{\sus}(Y;E))\cong
\bbZ.$ Using the connectedness of $G^{-\infty}(Y;E),$ we may
deform $f$ so that it is equal to the identity in a neighbourhood of some
point $p\in\bbS^{2k+1}$ and of the antipodal point $-p.$ This gives rise to
a new map $g:\bbR\times \bbS^{2k}\rightarrow G^{-\infty}(Y;E)$ with $g(s,
\omega)=\Id$ for $|s|>0$ large enough, where $(s,\omega)\in
\bbR\times\bbS^{2k}.$ By \eqref{fipomb.25}, the map $g$ can also be thought
as a map $g:\bbS^{2k}\rightarrow G^{-\infty}_{\sus}(Y;E).$ The argument can
be reversed, so this is a generator of
$\pi_{2k}(G^{-\infty}_{\sus}(Y;E))\cong \bbZ.$ It is now easy to relate
$c_k^{\even}$ to $c_k^{\odd}.$ If $(s,\omega)$ denotes coordinates for
$\bbR\times\bbS^{2k},$ we have
\begin{equation}
\begin{aligned}
1 &=\int_{\bbS^{2k+1}}f^{*}\beta_{k+1}^{\odd}=
\int_{\bbR\times\bbS^{2k}}g^{*}\beta_{k+1}^{\odd}
\\ 
&=\int_{\bbR\times\bbS^{2k}}g^{*}c_{k+1}^{\odd}
\Tr\left((a^{-1}d a)^{2k+1}\right)    \\
 &= c_{k+1}^{\odd}\int_{\bbR\times\bbS^{2k}}
\Tr\left((a^{-1}(s,\omega)\left(\frac{da}{ds}ds+
d_\omega a\right))^{2k+1}\right) \\
&=(2k+1)c_{k+1}^{\odd}\int_{\bbR\times\bbS^{2k}}
\Tr\left(\left(a^{-1}(s,\omega)d_\omega a\right)^{2k}
\left(a^{-1}(s,\omega)\frac{da}{ds}ds\right)\right) \\
&=(2k+1)\frac{c_{k+1}^{\odd}}{c_k^{\even}}
\int_{\bbS^{2k}} g^{*}\beta_k^{\even}=
(2k+1)\frac{c_{k+1}^{\odd}}{c_k^{\even}},  
\end{aligned}  
\end{equation}\label{fipombf.1}
so we conclude that $c_k^{\even}=(2k+1)c_{k+1}^{\odd}=\frac{1}{(2\pi
i)^{k+1}}\frac{k!}{(2k)!}$ .
\end{proof}

We define forms on $G^{-\infty}_{\cusp}(X;E)$ using the regularized trace
functional defined in \eqref{fipomb.31} by a choice of boundary defining
function and analytic continuation: 
\begin{equation}
\gamma_k=c_k^{\odd}\bTr\left((a^{-1}d a)^{2k-1}\right).
\label{fipomb.11}\end{equation}

\begin{proposition}\label{fipomb.12} If $X$ is a connected compact manifold
with boundary then for each $k\in\bbN,$  
\begin{equation}
d\gamma_k
= I^*(\beta_k^{\even}+dT_k),
\label{fipomb.14}\end{equation}
where $T_k$ is also a smooth form on $G^{-\infty}_{\sus}(\pa X;E)$ and 
\begin{equation}
I:G^{-\infty}_{\cusp}(X;E)\longrightarrow G^{-\infty}_{\sus}(\pa X;E)
\label{fipomb.20}\end{equation}
is the indicial map.
\end{proposition}

\begin{proof} Using \eqref{fipomb.10} and \eqref{fipomb.33}, an explicit
computation gives 
\begin{equation}
\begin{aligned}
d\gamma_k&=c_k^{\odd}\bTr\left(d(a^{-1}da)^{2k-1}\right) \\
&=-c_k^{\odd}\bTr\left((a^{-1}da)^{2k}\right) \\
&=-\frac12c_k^{\odd}\bTr\left([(a^{-1}da)^{2k-1},a^{-1}da]\right) \\
&=\frac12c_k^{\odd}\frac{1}{2\pi i}\int_{\bbR}
\Tr\left((a^{-1}da)^{2k-1}(a^{-1}\dot a
a^{-1}da-a^{-1}d\dot a)\right)d\tau\\ 
&=\frac{1}{4\pi i}c_k^{\odd}
\left(\int_{\bbR}\Tr\left((a^{-1}da)^{2k}a^{-1}\dot a\right)d\tau -
d\int_{\bbR}\Tr\left((a^{-1}da)^{2k-1}a^{-1}\dot a\right)d\tau\right) \\
&=I^*\left(\beta_k^{\even} 
+d\left(-\frac{1}{4\pi i}c_k^{\odd}\int_{\bbR}
\Tr\left((a^{-1}da)^{2k-1}a^{-1}\dot a\right)d\tau \right)\right).
\end{aligned}
\label{fipomb.13}\end{equation}
\end{proof}

Since $\dot\Psi^{-\infty}(X;E)\subset\Psi^{-\infty}_{\cusp}(X;E)$
is an ideal, $\dot G^{-\infty}(X;E)\subset G^{-\infty}_{\cusp}(X;E)$ is a
normal subgroup. The normal operator is trivial on $\dot\Psi^{-\infty}(X;E)$
and so maps the quotient into $G^{-\infty}_{\sus}(\pa X;E).$ Applying
Proposition~\ref{7.3.2003.7} shows that it has range exactly
$G^{-\infty}_{\sus,\ind=0}(\pa X;E).$ Thus we have a short exact sequence
of groups
\begin{equation}
\dot G^{-\infty}(X;E)\longrightarrow
G^{-\infty}_{\cusp}(X;E)\overset{N}\longrightarrow
  G^{-\infty}_{\sus,\ind=0}(\pa X;E)[[x]]
\label{5.3.2003.7}\end{equation}
where the quotient group is precisely the space of formal power series
\begin{multline}
G^{-\infty}_{\sus,\ind=0}(\pa X)[[x]]\\
=\{L\in\Psi^{-\infty}_{\sus}(\pa X;E)[[x]];
L=\sum\limits_{k=0}^\infty L_kx^k,\
\Id+L_0\in G^{-\infty}_{\sus,\ind=0}(\pa X;E)\}.
\label{fipomb.41}\end{multline}
The product here is a $\star$ product.

The lower order terms in the image group in \eqref{5.3.2003.7} are
arbitrary affine terms, so it is contractible to
$G^{-\infty}_{\sus,\ind=0}(X;E).$ The sequence \eqref{5.3.2003.7} is not
quite a fibration, since there is no local splitting. However, we have the
following weaker result.

\begin{lemma}\label{fipombf.3} The exact sequence of groups
\eqref{5.3.2003.7} is a Serre fibration, that is, it has the homotopy
lifting property for disks.
\end{lemma}

\begin{proof} Let $I=[0,1]$ be the unit interval and let
$h_{t}:I^{k}\longrightarrow G^{-\infty}_{\sus,\ind=0}(\pa X;E)[[x]],$ $t\in
[0,1],$ be a homotopy such that $h_{0}$ has a lift to
$G^{-\infty}_{\cusp}(X;E).$ Thus, there exists a map
$\tilde{h}_{0}:I^{k}\longrightarrow G^{-\infty}_{\cusp}(X;E)$ with the
property that $N\circ \tilde{h}_{0}=h_{0}.$ Then we must show that we can
extend $\tilde{h}_{0}$ to a homotopy $\tilde{h}_{t}:I^{k}\longrightarrow
G^{-\infty}_{\cusp}(X;E)$ such that $N\circ \tilde{h}_{t}=h_{t}$ for all
$t\in[0,1].$

First note that, using Borel's lemma, the homotopy $h_{t}$ can be lifted to
a homotopy $h_{t}':I^{k}\longrightarrow (\Id +\Psi^{-\infty}_{\cusp}(X;E))$
such that $N\circ h_{t}'=h_{t}$ for all $t\in[0,1].$ Initially however,
nothing ensures that the image of $h_{t}'$ lies in $G^{-\infty}_{\cusp}(X;E).$

Nevertheless, by Propositon~\ref{7.3.2003.7}, we know that the image of
$h_{t}'$ is contained in the space of Fredholm operators of index zero. In
fact, again by propositon~\ref{7.3.2003.7}, the homotopy $h_{t}'$ gives
rise to a principal bundle
\begin{equation}
\xymatrix{\dot
  G^{-\infty}(X;E)\ar@{-}\ar[r]&\mathcal{B}\ar[d]\\&[0,1]\times I^{k}} 
\label{fipombf.4}\end{equation}
 where the fibre $\mathcal{B}_{(t,x)}$ above the point $(t,x)\in [0,1]\times
 I^{k}$ is given by
\begin{equation}
\mathcal{B}_{(t,x)}= \{(h_{t}'(x) +A)\in G^{-\infty}_{\cusp}(X;E);
 A\in \dot\Psi^{-\infty}(X;E) \} . 
\label{fipombf.5}\end{equation}
Here, the group $\dot G^{-\infty}(X;E)$ has an obvious left action on the
fibre. To show that it has a principal bundle structure, we need to show that
there exist local trivializations everywhere on $[0,1]\times I^{k}$.   
For $(t_{0},x_{0})\in [0,1]\times I^{k}$ it follows from the vanishing of
the index of $h_{t_{0}}'(x_{0})$ that there exists 
$A\in\dot{\Psi}^{-\infty}(X;E)$
such that $h_{t_{0}}'(x_{0})+ A$ is invertible. Then $h_{t}(x)+A\in
G^{-\infty}_{cu}(X;E)$ for $(t,x)\in\mathcal{U},$ a small neighborhood of
$(t_{0},x_{0}).$ This provides a local trivialization of the bundle $\mathcal{B}$
\begin{equation}
\begin{matrix}
 F: &\mathcal{U}\times\dot{G}^{-\infty}(X;E)  &\longrightarrow 
 & \mathcal{B}  \\
            &  ((t,x),Q)     & \longmapsto    & (h_{t}(x)+A)\circ Q 
\end{matrix}
\end{equation}
and so \eqref{fipombf.4} has indeed a principal bundle structure.
 
Since the base space is contractible, this principal bundle must be
trivial (see for instance corollary 11.6 in \cite{MR2000a:55001}). In
particular, it has a global smooth section $s:[0,1]\times I^{k}\longrightarrow
\mathcal{B}.$  Without loss of generality, we can assume that this smooth
section is such that $s(0,x)=\tilde{h}_{0}(x)$ for all $x\in I^{k}.$
Finally, looking at the section $s$ as a map $s:[0,1]\times I^{k}\longrightarrow
G^{-\infty}_{\cusp}(X;E),$ we see that the homotopy $
\tilde{h}_{t}:I^{k}\longrightarrow G^{-\infty}_{\cusp}(X;E)$ defined by
$\tilde{h}_{t}(x)=s(t,x)$ is the desired lift of $h_{t}.$
\end{proof}

It follows that there is a long exact sequence in homotopy theory
\begin{multline}
\dots\longrightarrow 
\pi_{l}(\dot G^{-\infty}(X;E))
\longrightarrow
\pi_{l}(G^{-\infty}_{\cusp}(X;E))
\longrightarrow\\
\pi_{l}(G^{-\infty}_{\sus,\ind=0}(\pa X;E)[[x]])
\overset{b}\longrightarrow 
\pi_{l-1}(\dot G^{-\infty}(X;E))\dots ,
\label{fipomb.42}\end{multline} 
all classes being represented by smooth maps.

\begin{proposition}\label{5.3.2003.6} On any compact manifold, $X,$ with
non-trivial boundary, and for any complex vector bundle $E$ over
$X,$ the group $G^{-\infty}_{\cusp}(X;E)$ defined by \eqref{fipomb.40} is
weakly contractible.
\end{proposition}

\begin{proof} Using the long exact sequence
\eqref{fipomb.42}, we see that the weak contractibility of the group
$G^{-\infty}_{\cusp}(X;E))$ is equivalent to the fact that the maps
\begin{equation}
b:\pi_{l}(G^{-\infty}_{\sus,\ind=0}(\pa X;E)[[x]])\longrightarrow
\pi_{l-1}(\dot G^{-\infty}(X;E))
\label{fipombf.6}\end{equation}
are isomorphisms for all $l\in \bbN_{0}.$

For $l=0,$ this follows from the fact that $G^{-\infty}_{\sus,\ind=0}(\pa
X;E)[[x]]$ is connected (see theorem 4 in \cite{MR96h:58169}). For
$l=2k+1, k\in\bbN_{0},$ this is clearly the case, since
$\pi_{2k+1}(G^{-\infty}_{\sus,\ind=0}(\pa X;E)[[x]])= \pi_{2k}(\dot
G^{-\infty}(X;E))=\{0\}.$ Finally, for $l=2k, k\in\bbN,$
$\pi_{2k}(G^{-\infty}_{\sus,\ind=0}(\pa X;E)[[x]])\simeq \pi_{2k-1}(\dot
G^{-\infty}(X;E))\simeq \bbZ,$ so the map $b$ will be an isomorphism if it
maps a generator of $\pi_{2k}(G^{-\infty}_{\sus,\ind=0}(\pa X;E)[[x]])$ to
a generator of $\pi_{2k-1}(\dot G^{-\infty}(X;E)).$ Therefore, let us
assume without loss of generality that $f:\bbS^{2k}\longrightarrow
G^{-\infty}_{\sus,\ind=0}(\pa X;E))$ is such that $[f]\in
\pi_{2k}(G^{-\infty}_{\sus,\ind=0}(\pa X;E))\simeq
\pi_{2k}(G^{-\infty}_{\sus,\ind=0}(\pa X;E)[[x]])$ is a generator.

We proceed to construct $b([f]).$ First, lift $f$ to the 1-point
blow-up 
\begin{equation*}
\beta_{q}:[\bbS^{2k},\{q\}]\longrightarrow \bbS^{2k}.
\label{fipomb.65}\end{equation*}
Thus, $f\circ\beta_{q}$ is a
map from the $2k$-ball and consequently has a lift
$\tilde{f}:[\bbS^{2k},\{q\}]\longrightarrow G^{-\infty}_{\cusp}(X;E).$ Since
$I(\partial [\bbS^{2k};\{q\}])=q,$ $\tilde{f}$ has its image in a fibre
when restricted to $\partial [\bbS^{2k};\{q\}],$ so it is of the form
$\tilde{f}(s)=g(s)\circ \tilde{f}(q)$ for
$s\in\partial[\bbS^{2k};\{q\}]=\bbS^{2k-1}$ with $g:\bbS^{2k-1}\longrightarrow \dot
G^{-\infty}(X;E).$ This is just Serre's construction, see for instance
\cite{MR2000a:55001}, of $b([f])=[g]\in\pi_{2k-1}(\dot G^{-\infty}(X;E)).$
Hence, what we need to show is that $[g]$ is a generator of $\pi_{2k-1}(\dot
G^{-\infty}(X;E)).$  By Proposition~\ref{fipomb.6}, this amounts to
showing that $\int_{\bbS^{2k-1}}g^{*}\beta_{k}^{\odd}=\pm 1.$ By
Proposition~\ref{fipomb.7}, Proposition~\ref{fipomb.12} and Stokes'
theorem
\begin{equation}
\begin{aligned}
\int_{\bbS^{2k-1}}g^{*}\beta_{k}^{odd}&=&
\int_{\partial[\bbS^{2k},\{q\}]}\tilde{f}^{*}\gamma_{k}=
\int_{[\bbS^{2k},\{q\}]}\tilde{f}^{*}d\gamma_{k} \\ 
 &=&\int_{\bbS^{2k}}\tilde{f}^{*}(I^{*}(\beta_{k} +
dT_{k}))=\int_{\bbS^{2k}}f^{*}\beta_{k}=1 
\end{aligned}
\label{fipombf.2}\end{equation}
since $f$ was assumed to be a generator of
$\pi_{2k}(G^{-\infty}_{\sus,\ind=0}(\pa X;E)).$  Therefore, $b([f])=[g]$ is 
a generator of $\pi_{2k-1}(\dot G^{-\infty}(X;E)),$ which implies that the
map $b$ is an isomorphism and we conclude that the group
$G^{-\infty}_{\cusp}(X;E)$ is weakly contractible.
\end{proof}

\section{Perturbations of elliptic cusp operators}

\begin{proposition}\label{7.3.2003.8} For any elliptic operator
  $P\in\Psi_{\cusp}^{m}(X;E,F)$ there exist perturbations
  $Q\in\Psi_{\cusp}^{-\infty}(X;E,F)$ such that $P+Q$ is Fredholm as a
  map from $\dCI(X;E)$ to $\dCI(X;F),$ or from the cusp Sobolev space
  $H^m_{\cusp}(X;E)$ to $L^2_{\cusp}(X;F).$
\end{proposition}

\begin{proof} From the definition of
  $\Psi_{\cusp}^{m}(X;E)$ we may cut the kernel of $P$ off near the boundary of
  the front face corresponding to each pair of boundary components and so arrange
  that the normal operator has compact support. This only changes $P$ by an
  element of $\Psi_{\cusp}^{-\infty}(X;E)$ and ensures that the indicial family
  is an entire family of elliptic elements $I(P,s)\in\Psi^{m}(\pa X;E).$ In
  fact this family is necessarily invertible as $|\Re s|\to\infty$ with $\Im
  s$ bounded. We may perturb it to make it invertible on the
  real axis, which implies that the resulting operator is
  Fredholm. Indeed, $I(P,s)$ is an entire, Fredholm family which is
  invertible at one point, so the set of points at which it has non-trivial
kernel is discrete. Near such a point $s_0\in\bbR,$ we may divide the
  space, either $\dCI(\pa X;E)$ or $H^m(\pa X;E),$ into the finite
  dimensional null space and its $L^2$ orthocomplement. With a similar
  decomposition of the range in terms of the adjoint, $I(P,s)$ becomes a
  matrix of the form
\begin{equation}
\begin{pmatrix}P'&(s-s_0)Q\\(s-s_0)L&(s-s_0)A \end{pmatrix}.
\label{7.3.2003.9}\end{equation}
The matrix $M(s)=\begin{pmatrix}P'&0\\ 0&G\end{pmatrix},$ where $G$ is an
isomorphism from the null space to the null space of the adjoint, is
invertible. Composing on the left with the inverse gives a matrix 
\begin{equation}
\begin{pmatrix}\Id&(s-s_0)Q'\\(s-s_0)L'&(s-s_0)A' \end{pmatrix}.
\label{7.3.2003.10}\end{equation}
This matrix is of the form $\Id+R(s)$ with $R(s)$ of finite rank, and
smoothing. Moreover the determinant vanishes only to finite order at
$s=s_0$ so there is a smoothing perturbation supported arbitrarily
close to $s_0$ making it invertible in a fixed neighbourhood of $s_0.$
For instance, choosing $\epsilon>0$ small enough so that 
$\Id +R(s)$ is invertible in the punctured disc $D_{\epsilon}(s_{0})-\{s_{0}\}$, where
\begin{equation}
D_{\epsilon}(s_{0})=\{ s\in \bbC ; |s-s_{0}|<\epsilon \}\, ,
\label{disc}\end{equation}
the perturbation $\Id+ R(s+i\varphi(s))$ will do, where 
$\varphi\in \CI_{c}(\bbR)$
is a real-valued function satisfying
\begin{equation}
  -\frac{\epsilon}{2}<\varphi(s)<\frac{\epsilon}{2}\
\forall\ s\in \bbR, \; \; \varphi(s_{0})\ne 0\text{ and }
\supp{\varphi}\subset [s_{0}-\frac{\epsilon}{2},s_{0}+\frac{\epsilon}{2}].
 \label{pertubation}\end{equation}

Doing this at each point at which there is non-trivial null space
gives a perturbation of $\Id+R(s)$ making it invertible. Composing on
the left with $M(s)$ gives a perturbation of $I(P,s)$ making it
invertible for all $s\in\bbR.$ This perturbation is the indicial family of an
element $Q\in\Psi_{\cusp}^{-\infty}(X;E,F)$ which therefore gives the
desired perturbation of $P.$
\end{proof}

We may strengthen Proposition~\ref{7.3.2003.8} further.

\begin{corollary}\label{CnvRcn.6} Under the assumptions of
  Proposition~\ref{7.3.2003.8} there exists a perturbation
  $Q\in\Psi_{\cusp}^{-\infty}(X;E,F)$ such that $P+Q$ is an
  isomorphism, in either sense.
\end{corollary}

\begin{proof} Having perturbed $P$ by a term 
  $Q_1\in\Psi_{\cusp}^{-\infty}(X;E,F)$ so that it is Fredholm,
  suppose its index is $k\in\bbZ.$ Using Proposition~\ref{7.3.2003.7}, we
  may multiply on the right by an operator of the form $\Id+A,$
  $A\in\Psi_{\cusp}^{-\infty}(X;E)$ with index $-k.$ This gives an operator
  $P+Q_2,$ still with $Q_2\in\Psi_{\cusp}^{-\infty}(X;E,F),$ which is
  Fredholm of index zero. The null space is in $\dCI(X;E)$ and the null
  space of the adjoint, with respect to a choice of fibre metric, is in
  $\dCI(X;F)$ and has the same dimension. Thus, adding a finite rank smoothing
  operator in $\dot\Psi(X;E,F)$ gives an invertible operator $P+Q$ with $Q$
  as claimed.
\end{proof}

\section{Perturbation of Elliptic Families}

Now consider a smooth fibration of compact manifolds \eqref{fipomb.1}
where the model fibre $Z$ is a manifold with boundary. If $P\in
\Psi_{\cusp}^m(M/B;E,F)$ is a smooth family of elliptic cusp operators
on the fibres acting between sections of vector bundles over the
total space $M,$ then for $b\in B$ we may consider 
\begin{equation}
\mathcal{P}_b=\{P_b+Q_b; Q_b\in\Psi_{\cusp}^{-\infty}(Z_b;E_b,F_b),\ \exists\
(P_b+Q_b)^{-1}\in \Psi_{\cusp}^{-m}(Z_b;F_b,E_b)\}.
\label{CnvRcn.9}\end{equation}
Here $Z_b=\phi^{-1}(b)$ is the fibre and $E_b,$ $F_b$ are the
restrictions of the bundles to it.

\begin{proposition}\label{CnvRcn.8} For an elliptic family $P\in
\Psi_{\cusp}^m(M/B;E,F)$ the spaces in \eqref{CnvRcn.9} form a smooth
bundle of principal $G^{-\infty}_{\cusp}(Z_b;E)$-spaces over $B,$ 
\begin{equation}
\xymatrix{G^{-\infty}_{\cusp}(Z_b;E)\ar@{-}[r]&\mathcal{P}\ar[d]^{\phi}\\
&B.}
\label{CnvRcn.10}\end{equation}
\end{proposition}

\begin{proof} By Corollary~\ref{CnvRcn.6} the spaces in
  \eqref{CnvRcn.9} are non-empty. Moreover, the notion of a local smooth
  section is well defined. Furthermore, any two points in
  $\mathcal{P}$ corresponding to invertible perturbations by $Q_1$
  and $Q_2$ are related by composition on the right with
  $(P_b+Q_2)^{-1}(P_b+Q_1)\in G^{-\infty}_{\cusp}(Z_b;E_b)$ and conversely.
\end{proof}

Using the contractibility of the fibres, it is now easy to show the
existence of a global section of the bundle $\mathcal{P}.$

\begin{theorem}\label{5.3.2003.1} Let $P$ be a smooth family of elliptic
cusp pseudodifferential operators of order $m$ on the fibres of a fibration
of a compact manifold with boundary. Then there exists a smooth family $Q$
of cusp pseudodifferential operators of order $-\infty$ such that
$P+Q:\dCI(X;E)\longrightarrow\dCI(X;F)$ is invertible.
\end{theorem}

\begin{proof}
Given $p\in B,$ there exists an open neighbourhood $\mathcal{U}$ of $p$ and
a local trivialization of the bundle $M,$ 
\begin{equation}
h:\mathcal{U}\times Z \longrightarrow \phi^{-1}(\mathcal{U}),\ \ with\ \
\phi\circ h(u,z)=u \ \forall\ (u,z)\in \mathcal{U}\times Z.
\label{CnvRcn.R11}\end{equation}
This gives rise to a local trivialization of the bundle $\mathcal{P},$ and
in fact to a local structure of principal $G^{-\infty}_{\cusp}(Z;E)$-bundle
\begin{equation}
H:\mathcal{U}\times G^{-\infty}_{\cusp}(Z;E) \longrightarrow
\phi^{-1}(\mathcal{U})\ .
\label{CnvRcn.R12}\end{equation}

Notice however that this local structure of principal
$G^{-\infty}_{\cusp}(Z;E)$-bundle depends on the choice of the
trivialization in \eqref{CnvRcn.R11}, so is not canonical. In particular,
this means that there is no natural way of defining a global
principal bundle structure with fixed structure group. Nevertheless this
suffices to construct a global continuous section of the bundle $\mathcal{P}.$
Indeed, consider a triangulation of the manifold $B$ such that each simplex
is contained in some open set $\mathcal{U}$ as in \eqref{CnvRcn.R12}. Then,
assigning arbitrary values on the vertices of this triangulation to the
section $s:B\longrightarrow \mathcal{P}$ we want to construct, we see,
using Proposition~\ref{5.3.2003.6} and the discussion of $\S 29.1$ in
\cite{MR2000a:55001}, that it is possible to extend the section $s$
continuously to all of $B.$ Using a partition of unity and local smoothing
we may perturb this continuous section to be smooth.
\end{proof}

\section{Relative index theorem}

Under the same, ellipticity, hypothesis as in Theorem~\ref{5.3.2003.1} we
now consider all the Fredholm perturbations in the cusp algebra.

\begin{theorem}\label{5.3.2003.3} If $P\in\Psi_{\cusp}^{m}(M/B;E,F)$
is a smooth family of elliptic cusp pseudodifferential operators of order
$m$ on the fibres of a fibration of a compact manifold with boundary then
\begin{equation}
\begin{gathered}
\pi_0(\mathcal{P}_F)\simeq K^0(B)\\
\mathcal{P}_F=
\left\{Q\in\Psi_{\cusp}^{-\infty}(M/B;E,F);P+Q
\text{ is a family of Fredholm operators}\right\}
\end{gathered}
\label{fipomb.46}\end{equation}
with the isomorphism given by the relative index formula
\begin{equation}
\ind(P+Q)=[I(P+Q)I(P+Q_0)^{-1}]\Min K^0(B)
\label{5.3.2003.2}\end{equation}
where $P+Q_0$ is invertible.
\end{theorem}

\begin{proof} We have shown above that the bundle $\mathcal{P}$ in
\eqref{CnvRcn.10} has a global section consisting of invertible operators;
$P+Q_0$ in \eqref{5.3.2003.2} is such a choice. The space of Fredholm
perturbations, which appears in \eqref{fipomb.46} consists precisely of
those $Q$ such that $I(P+Q)\in\Psi_{\sus}^m(\pa M/B;E,F)$ is
invertible. Since $I(P+Q_0)$ is invertible this defines a map 
\begin{equation}
\mathcal{P}_F\ni Q\longmapsto I(P+Q)I(P+Q_0)^{-1}\in G^{-\infty}_{\sus}(\pa
M/B;F).
\label{fipomb.47}\end{equation}
Moreover the fibres of this map are contractible, since if $I(P+Q')=I(P+Q)$
with $Q,Q'\in\mathcal{P}_F$ then $sQ'+(1-s)Q\in\mathcal{P}_F.$ Thus the
families index only depends on this normalised indicial
family. Furthermore, \eqref{fipomb.47} is surjective so 
\begin{equation}
\pi_0(\mathcal{P}_F)\simeq \pi_0(G^{-\infty}_{\sus}(\pa M/B;F))\simeq K^0(B).
\label{fipomb.48}\end{equation}
\end{proof}

\section{Odd case}

The suspended algebra of cusp operators can be defined by close analogy
with the suspension of the usual algebra. That is, for a compact manifold
with boundary $Z$ and vector bundles $E,$ $F,$ consider the space of cusp
operators $\Psi^m_{\cusp}(\bbR\times Z;E,F)$ on $\bbR\times Z$
corresponding to the lift of a cusp structure from $Z.$ Thus the algebra is
translation-invariant and within it consider the subspace of
translation-invariant operators. Away from the diagonal the kernels of
these operators are smooth functions on 
\begin{equation}
(\bbR\times Z)^2_{\cusp}=\bbR\times\bbR\times Z^2_{\cusp}.
\label{fipomb.61}\end{equation}
Thus we may further consider the subspace of kernels which are Schwartz
near infinity. This subspace is closed under composition and we denote the
resulting filtered algebra $\Psi^m_{\scusp}(Z;E,F).$ 

The symbol sequence for these operators  
\begin{equation}
\Psi^{m-1}_{\scusp}(Z;E,F)\longrightarrow \Psi^m_{\scusp}(Z;E,F)\longrightarrow
\CI(S(\bbR\times T^*_{\cusp}Z);\hom(E,F)\otimes R^m)
\label{fipomb.62}\end{equation}
has the usual multiplicative properties as does the boundary sequence 
\begin{equation}
x\Psi^m_{\scusp}(Z;E,F)\longrightarrow \Psi^m_{\scusp}(Z;E,F)
\longleftrightarrow \Psi_{\sus(\bbR\times L)}^m(\pa Z;E,F)
\label{fipomb.63}\end{equation}
where $\Psi_{\sus(\bbR\times L)}^m(\pa Z;E,F)$ is the doubly-suspended algebra of
pseudodifferential operators on the boundary.

The `odd' families index theorem can be viewed as the families index
theorem for this algebra. In essence this is Theorem~\ref{5.3.2003.1}
except that we are dealing with a special family.

\begin{theorem}\label{fipomb.21} For any elliptic suspended family on the
fibres of a fibration of a compact manifold with boundary,
$P\in\Psi^m_{\scusp}(M/B;E),$ there is a perturbation
$Q\in\Psi^{-\infty}_{\scusp}(M/B;E)$ such that $P+Q$ is invertible with
inverse in $\Psi^{-m}_{\scusp}(M/B;E)$ and the set of homotopy classes of
Fredholm perturbations $Q\in\Psi^{-\infty}_{\scusp}(M/B;E)$ is naturally
isomorphic to $K^1(B).$
\end{theorem}

\begin{proof} Taking the Fourier transform in the suspension variable gives
a 1-parameter family of suspended operators. If the family is elliptic this
is automatically invertible for large enough values of the
parameter. Restricting the parameter space to $[-R,R]\times B$ for $R$
large enough, we may apply Proposition~\ref{7.3.2003.8} to find a
perturbation $Q'\in\Psi_{\scusp}^{-\infty}(M/B;E,F)$ making it invertible
over $[-R,R]\times B.$ Restricting the parameter to values $\pm R$ gives,
in each case, two invertible families over $B,$ namely $P+Q$ and $P.$ Since 
\begin{equation}
(P+Q)P^{-1}=\Id+R,\ R\in \Psi_{\cusp}^{-\infty}(\pa M/B;F)
\label{fipomb.64}\end{equation}
and we know weak contractibility of this group from
Proposition~\ref{5.3.2003.6}, we may modify the family of perturbations $Q'$
to vanish outside $[-R-1,R+1]\times B$ and still have $P+Q'$ everywhere
invertible. Taking the inverse Fourier transform in the suspending variable
gives the invertible perturbation.

If we consider perturbations making the family everywhere Fredholm, \ie so
that the suspended indicial family (so doubly-suspended families over the
boundary) is everywhere invertible, we get an odd index, with values in
$K^1(B).$ Two different perturbations of this type correspond to a family
of double-suspended invertible operators in $G^{-\infty}(\pa Z)$ and hence
to an element of $K^1(B)$ by Bott periodicity. This is our relative index
theorem, which shows that the set of components of the space of everywhere
Fredholm perturbations is isomorphic to $K^1(B).$
\end{proof}

\appendix

\section*{Appendix:Cusp algebra}

On a compact manifold with boundary there are various distinct classes of
pseudodifferential operators. Here we use the algebra of `cusp'
operators. These operators can be thought of as `asymptotically
translation-invariant' when a neighbourhood of the boundary is identified,
through inversion of a defining function, with a `cylindrical end'. The
definition and development of the basic properties of these operators is
due to Rafe Mazzeo and the first author; more details (and a more general
class of operators) can be found in \cite{Mazzeo-Melrose4}. Notice however,
that the definition here is slightly different to that in
\cite{Mazzeo-Melrose4} in that we use the `overblown' cusp calculus, which
admits terms relating the different boundary components.

First consider the appropriate class of translation-invariant
operators. Thus, if $Y$ is a compact manifold without boundary we may
consider pseudodifferential operators on $\bbR\times Y,$ acting from
sections of one vector bundle, $E$ over $Y$ to another, $F.$ If
$A\in\Psi^m(\bbR\times Y;E,F)$ is such an operator which is invariant under
translation in $\bbR$ then it may be considered as a convolution operator
on $\bbR$ and hence has as Schwartz kernel a distribution,
$A\in\CmI(\bbR\times Y^2;\Hom(E,F)\otimes\pi^*_R\Omega(Y))$ which is smooth 
away from the `diagonal' $\{0\}\times\Diag_Y$ and has a classical conormal
singularity at the diagonal. As in \cite{MR96h:58169}, we consider the
subspace of those operators which have kernels which are rapidly decreasing
at infinity with all derivatives. These operators, the space of which
we denote $\Psi^m_{\sus}(Y;E,F),$ map $\mathcal{S}(\bbR\times Y;E)$ into
$\mathcal{S}(\bbR\times Y;F).$ They have the natural composition property 
\begin{equation}
\Psi^m_{\sus}(Y;F,G)\circ \Psi^{m'}_{\sus}(Y;E,F)=\Psi^{m+m'}_{\sus}(Y;E,G).
\label{fipomb.22}\end{equation}
Fourier transformation on $\bbR$ converts $A\in\Psi^m_{\sus}(Y;E,F)$ into a
family of pseudodifferential operators on $Y:$ 
\begin{equation}
A(e^{it\tau}u)=e^{it\tau }\hat A(\tau)u,\ u\in\CI(Y;E).
\label{fipomb.25}\end{equation}
The family $\hat A:\bbR\longrightarrow \Psi^m(Y;E,F)$ determines $A$ but is
not an arbitrary family of pseudodifferential operators on $Y$ depending
smoothly on $\tau\in\bbR.$ Rather, the families arising this way may be
characterised as having full symbols, in terms of some global quantisation,
which are (classical) symbols on $\bbR_\tau \times T^*Y.$

If $Y_1$ and $Y_2$ are different manifolds then the discussion above
extends directly to define the spaces $\Psi^{-\infty}_{\sus}(Y_1,Y_2;E_1,E_2)$
of suspended smoothing operators between sections of bundles $E_i$ over $Y_i,$
$i=1,2.$ In general there is no analogous space of positive order
pseudodifferential operators. These smoothing operators form modules over
the spaces of suspended pseudodifferential operators: 
\begin{equation}
\begin{gathered}
\Psi^{-\infty}_{\sus}(Y_1,Y_2;E_1,E_2)\circ \Psi^{m}_{\sus}(Y_1;E_0,E_1)
\subset\Psi^{-\infty}_{\sus}(Y_1,Y_2;E_0,E_2),\\
\Psi^{m}_{\sus}(Y_2;E_2,E_3)\circ\Psi^{-\infty}_{\sus}(Y_1,Y_2;E_1,E_2)
\subset\Psi^{-\infty}_{\sus}(Y_1,Y_2;E_1,E_3),\\
\Psi^{-\infty}_{\sus}(Y_2,Y_3;E_2,E_3)\circ
\Psi^{-\infty}_{\sus}(Y_1,Y_2;E_1,E_2)\subset
\Psi^{-\infty}_{\sus}(Y_1,Y_3;E_1,E_3).
\end{gathered}
\label{fipomb.67}\end{equation}
If $Y$ is a `disconnected manifold', meaning the union of a finite number
of disjoint smooth manifolds (of fixed dimension), then we will denote by
$\Psi^m_{\sus}(Y;E,F)$ the space of matrices, with entries parameterized by
$\pi_0(Y)$ of the type discussed above between the components of $Y.$ The
entries between different components are necessarily smoothing operators.

Now, let $X$ be a compact manifold with boundary, and suppose $E$ and $F$
are vector bundles over $X.$ We will assume initially that the boundary of
$X$ is connected. The most basic class of operators we use below, which we
denote by $\dot\Psi^{-\infty}(X;E,F),$ consists of those operators which
have Schwartz kernels in $\dCI(X^2;\Hom(E,F)\otimes\pi^*_R\Omega (X)),$
\ie these are smoothing operators with kernels vanishing to infinite order at
all boundary points of $X^2.$ They form an algebra which is isomorphic to
the algebra of smoothing operators on any manifold without boundary
(provided it is non-trivial, connected and of positive dimension).

If $0\le x\in\CI(X)$ is a boundary defining function and
$F:\{x<\epsilon\}\longrightarrow[0,\epsilon)\times\pa X$ is an associated
product decomposition near the boundary, for some $\epsilon >0,$ then we
may allow elements of $\Psi_{\sus}^m(\pa X;E,F)$ to act on $X$ by choosing a
cut-off function $\chi\in\CIc([0,\epsilon )),$ $\chi (x)=1$ in
$x\le\frac\epsilon 2,$ and setting
\begin{multline}
A'u=F^*(\chi (x)(Av)(1/x,\cdot)),\ v(t,\cdot)=\chi
(1/t)((F^{-1})^*u)(1/t,\cdot),\\ 
A\in\Psi_{\sus}^m(\pa X;E,F),\ u\in\dCI(X;E).
\label{fipomb.24}\end{multline}
Here it should be observed that $v\in\mathcal{S}(\bbR_t\times\pa X).$ This
action is not multiplicative.

Now the cusp pseudodifferential operators $\Psi_{\cusp}^m(X;E,F)\subset 
\Psi^m(\inn X;E,F)$ are operators $B:\CI(X;E)\longrightarrow
\CI(X;F)$ such that $gBg'\in\dot\Psi^{-\infty}(X;E,F),$ if $g,g'\in\CI(X)$
have supports with disjoint interiors and such that there exist elements
$A_j\in\Psi_{\sus}^m(\pa X;E,F)$ for which
\begin{multline}
B_\tau=e^{-i\tau /x}B e^{i\tau /x}:\CI(X;E)\longrightarrow \CI(X;F),\\
B_\tau \chi (x)F^*(u)\sim\sum\limits_{j=0}^\infty x^jF^*(\hat A_j(\tau )u)\
\forall\ u\in\CI(Y;E),\ \tau \in\bbR.
\label{fipomb.26}\end{multline}
A choice of defining function, $x,$ is involved here and, as noted
above, the algebra does depend on this choice, although only through the
trivialization of the normal bundle it determines. The algebra may be
characterised by an appropriately uniform version of \eqref{fipomb.26},
although this is a rather cumbersome approach.

Instead the Schwartz kernels of its elements may be characterised directly
on a cusp configuration space
\begin{equation}
X^2_{\cusp} =[X^2;\pa X\times\pa X;D_{\cusp}],\
\beta_{\cusp}:X^2_{\cusp}\longrightarrow X^2,
\label{fipomb.50}\end{equation}
which is a blown-up version of the usual product $X^2.$ Here $\pa
X\times\pa X$ is the corner of $X^2.$ If $\pa X$ is not connected it should
be replaced by the finitely many products between pairs of boundary
hypersurfaces. Note that the description here diverges slightly from that
in \cite{Mazzeo-Melrose4} where only the products of boundary hypersurfaces of
$X$ with themselves are considered and off-diagonal terms are not blown
up. The submanifold $D_{\cusp}\subset\ff([X^2;\pa X\times\pa X])$ is
determined by the cusp structure, in case $\pa X$ is not connected there is
one component corresponding to each pair of boundary faces. Thus, let $x_H$ be an
admissible defining function (corresponding to the chosen cusp structure)
for a boundary face $H.$ Then if $x'_{G}$ is the same function in a second
copy of $X,$ for a possibly different boundary hypersurface, the function 
\begin{equation}
s=\frac{x-x'}{x+x'}\in\CI([X^2;\pa X\times\pa X]),\ x=x_H,\ x'=x'_{G}
\label{fipomb.51}\end{equation}
has the property that $s=\pm1$ is the lift to $[X^2;\pa X\times\pa X]$ of
the two local boundary faces, $H$ as a boundary face in the first factor of
$X^2$ and $G$ as a boundary face in the second factor and, furthermore,
$\rho=\frac12(x+x')$ defines the new boundary hypersurface introduced in
the blow up. In fact $[X^2;\pa X\times\pa X]$ is precisely the space
obtained by appending $s$ to $\CI(X^2)$ (for each pair if boundary
hypersurfaces). Then
\begin{equation}
D _{\cusp}=\{s=0,\rho =0\}\subset\ff([X^2;\pa X\times\pa X])
\label{fipomb.52}\end{equation}
is determined by (and determines) the cusp structure.

The configuration space $X^2_{\cusp}$ has three type of boundary
hypersurfaces. The `old' boundary hypersurfaces, in two-to-one
correspondence with the boundary faces of $X,$ and three hypersurfaces
(because the second blow-up disconnects the hypersurfaces in the first blow
up), for each pair of boundary faces of $X,$ corresponding to the two blow-ups in
\eqref{fipomb.51}. We denote the last of these $\ff_{\cusp}.$ 
The diagonal terms (corresponding to a product of a boundary hypersurface
with itself) in this set are the only boundary hypersurfaces which meet the
diagonal $\Diag_{\cusp}\subset X^2_{\cusp},$ defined as the closure of the
inverse image of the interior of the diagonal of $X^2,$ and they meet
transversally. The part of the boundary other than $\ff_{\cusp}$ will be
denoted $B_{\nd}$ ($\nd$ for non-diagonal). Notice that the interior of the
cusp front face is of the form
\begin{equation}
\inn\ff_{\cusp}\equiv\bbR_t\times\pa X\times\pa X
\label{fipomb.59}\end{equation}
where the function $t=(x/x'-1)/\rho$ gives a coordinate on the first factor. 
The kernels of the residual part of the cusp algebra are by definition
those smooth densities on the interior of $X^2$ which are of the form
\begin{equation}
\Psi^{-\infty}_{\cusp}(X)=\left\{A=A'\nu_{\cusp};A'\in\CI(X^2_{\cusp}),\
A'\equiv0\Mat B_{\nd}\right\}.
\label{fipomb.54}\end{equation}
Here $\nu_{\cusp}$ is a smooth cusp density on the right fact of $X,$ \ie
of the form $\nu /x^2$ where $\nu$ is a smooth density in the usual sense.

This space of residual kernels is a module over $\CI(X^2_{\cusp})$ and can
be thought of as the `coefficient module'. This allows us to define related
spaces of kernels directly. Most importantly, for any $m\in\bbR$ we set 
\begin{equation}
\Psi^{m}_{\cusp}(X)=I^m(X^2_{\cusp};\Diag_{\cusp})\otimes_{\CI(X^2_{\cusp})}
\Psi^{-\infty}_{\cusp}(X)
\label{fipomb.55}\end{equation}
using the space of (classical) conormal distributions at the lifted
diagonal. The expansion \eqref{fipomb.26} corresponds precisely to the
expansion in Taylor series of the kernels at $\ff_{\cusp}.$ In particular
the leading term is determined by, and determines, the restriction of $A'$
to $\ff_{\cusp}.$

The leading term of the expansion in \eqref{fipomb.26} is independent of
choices, up to a positive scaling in the translation variable
(constant on the factor $\pa X).$ Thus, with only a choice of linear
variable on an oriented one-dimensional vector space, we may identify the
normal operator as a surjective algebra homomorphism
\begin{equation}
\Psi_{\cusp}^m(X;E,F)\longrightarrow \Psi_{\sus}^m(\pa X;E,F),\ B\longmapsto A_0
\label{fipomb.27}\end{equation}
defined through \eqref{fipomb.26}. We call the associated family of
pseudodifferential operators $I(B,\tau)=\hat A_0(\tau )\in\Psi(\pa Y;E,F)$ the
\emph{indicial family}. Note that when $\pa X$ is not connected the image
space here consists of the `matrices' labelled by two sets of the boundary
components of $\pa X,$ with entries which are suspended pseudodifferential
operators between the components as discussed above. When the components
are different the operators are necessarily smoothing, despite the formal
presence of the order indicator $m.$ 

As noted above, the choice of a cusp structure on a compact manifold with
boundary fixes a cusp tangent and cotangent bundle, isomorphic to the usual
tangent and cotangent bundles, but not naturally so. Namely, if $x$ is an
admissible defining function for the boundary then the smooth vector
fields, $V$ on $X$ satisfying 
\begin{equation}
Vx\in x^2\CI(X)\Longleftrightarrow V\in\mathcal{V}_{\cusp}(X)
\label{fipomb.34}\end{equation}
form a Lie algebra and $\CI(X)$-module depending only on the underlying
cusp structure and determining it. In local coordinates $x,$ $y_i,$ near a
boundary point, $\mathcal{V}_{\cusp}(X)$ is spanned over $\CI(X)$ by
$x^2\pa_x,$ $\pa_{y_i}.$ As a locally free module, this defines a smooth
vector bundle over $X$ (including over the boundary) such that 
\begin{equation}
\mathcal{V}_{\cusp}(X)=\CI(X;{}^{\cusp}TX).
\label{fipomb.35}\end{equation}

The dual bundle, ${}^{\cusp}T^*X,$ is the natural carrier for the symbols
of cusp pseudodifferential (or differential) operators. In particular, the
standard symbol map over the interior of $X$ extends by continuity to define
\eqref{CnvRcn.22} over any compact manifold with boundary. There is always
a global quantisation map 
\begin{equation}
q:S^m_{\cl}({}^{\cusp}T^*X)\longrightarrow \Psi^m_{\cusp}(X)
\label{fipomb.36}\end{equation}
(and similarly for operators between sections of vector bundles) which
induces a `full' symbol isomorphism 
\begin{equation}
S^m_{\cl}({}^{\cusp}T^*X)/S^{-\infty}({}^{\cusp}T^*X)\simeq
\Psi^m_{\cusp}(X)/\Psi^{-\infty}_{\cusp}(X) .
\label{fipomb.37}\end{equation}

The expansion \eqref{fipomb.26} is not trivial even in the case of
operators of symbolic order $-\infty.$ In fact, it then gives a short exact
sequence
\begin{equation}
\dot\Psi^{-\infty}(X;E,F)\longrightarrow\Psi^{-\infty}_{\cusp}(X;E,F)
\longrightarrow\Psi^{-\infty}_{\sus}(\pa X;E,F)[[x]].
\label{fipomb.17}\end{equation}
Here, the image consists of countably many copies of
$\Psi^{-\infty}_{\sus}(\pa X;E,F)$ with an induced `star' product
(if $E=F$ or involving three bundles) which is only trivially
multiplicative at the top level.

The elements of $\dot\Psi^{-\infty}(X;E)$ are of trace class but in general
those of $\Psi^{-\infty}_{\cusp}(X;E)$ are not, precisely because of the
cusp density in \eqref{fipomb.54}, which is not integrable. For $z\in\bbC$ and
$A\in\Psi^{-\infty}_{\cusp}(X;E),$
$x^zAx^{-z}=B_z\in\Psi^{-\infty}_{\cusp}(X;E).$ Moreover, for $\Re z>1,$
$x^zA$ is then of trace class and the function 
\begin{equation}
f_A(z)=\Tr(x^zA)=\int_{\Diag_{\cusp}} x^zA|{\Diag_{\cusp}}
=\int_{X} x^zA'|{\Diag_{\cusp}}\nu_{\cusp},\ A\in\Psi^{-\infty}_{\cusp}(X;E)
\label{fipomb.28}\end{equation}
is holomorphic in $\Re z>1$ with a meromorphic extension to the complex
plane. Since the cusp density is of the form $\nu /x^2$ in terms of a
smooth density, the possible poles are at the points $1-\bbN_0,$
$\bbN_0=\{0,1,\dots,\}$ and are at most simple. The \emph{boundary residue
trace} is defined as the residue at $z=0:$ 
\begin{equation}
\RTr\pa(A)=\lim_{z\longrightarrow0}zf_A(z):\Psi^{-\infty}_{\cusp}(X;E)\longrightarrow \bbC.
\label{fipomb.29}\end{equation}
It is independent of the choice of admissible defining function used to
define it, is a trace and is given explicitly in terms of the second term
in the expansion \eqref{fipomb.26}
\begin{equation}
\begin{gathered}
\RTr\pa([A,B])=0\ \forall\ A,B\in \Psi^{-\infty}_{\cusp}(X;E),\\
\RTr\pa(A)=\frac1{2\pi}\int_{\bbR}\Tr(\hat A_{1}(\tau))d\tau.
\end{gathered}
\label{fipomb.30}\end{equation}

The regularized value of $f_A(z)$ at $z=0$ is the \emph{regularized trace} 
\begin{equation}
\bTr(A)=\lim_{z\to0}\left(f_A(z)-\frac1z\RTr\pa(A)\right).
\label{fipomb.31}\end{equation}
It does depend on the choice of admissible defining function $x$ and is not
a trace, although it reduces to the usual trace on the trace class elements 
\begin{equation}
\bTr(A)=\Tr(A)\ \forall\ A\in x^2\Psi^{-\infty}_{\cusp}(X;E).
\label{fipomb.32}\end{equation}
In this paper, we make substantial use of the \emph{trace-defect formula}.

\begin{lemma}\label{fipomb.56} The regularized trace of a commutator only
depends on the indicial families (or the normal operators)
\begin{equation}
\bTr([A,B])=\frac1{2\pi i}\int_{\bbR}\Tr
\left(I(A,\tau)\frac{\pa}{\pa\tau}I(B,\tau)\right)d\tau,\
A,B\in\Psi^{-\infty}_{\cusp}(X;E).
\label{fipomb.33}\end{equation}
\end{lemma}

\begin{proof} This is essentially the same as the corresponding formula in
the b-calculus and can also be found in \cite{Melrose-Nistor2}. From the
definition, for $\Re z>>0,$
\begin{multline}
\Tr(x^z[A,B])=\Tr([x^z,A]B)=\Tr(x^z A_zB),\ A_z=x^{-z}[x^z,A]
\Longrightarrow\\
\bTr([A,B])=-\RTr\pa((D_{\log x}A)B).
\label{fipomb.57}\end{multline}
Here, $A_z\in\Psi^{-\infty}_{\cusp}(X;E)$ is entire in $z$ and vanishes at
$z=0$ and by definition 
\begin{equation}
D_{\log x}A=-\frac{\pa}{\pa z}A_z\big|_{z=0}.
\label{fipomb.58}\end{equation}
This is readily computed; the Schwartz kernel of $A_z$ is $A(1-(x'/x)^z)$
so the $z$-derivative at $z=0$ has kernel $A\log(x/x').$ Notice that $s=0$
in \eqref{fipomb.52} is $x/x'=1,$ $x+x'=0.$ Thus $\log(x/x')$ is a smooth
function, away from $x=0$ or $x'=0,$ which vanishes on $D_{\cusp}.$ Lifted
to $X^2_{\cusp}$ it is therefore of the form $\rho t$ where $t$ is given
following \eqref{fipomb.59}. Thus in fact 
\begin{equation}
D_{\log x}A\in x\Psi^{-\infty}_{\cusp}(X;E),\ I_{1}(D_{\log x}A)=
-i\frac{d}{ds}I(A,s)
\label{fipomb.60}\end{equation}
where $I_{1}(D_{\log x}A)=\widehat{(D_{\log{x}}A)}_{1}$ is the second term 
in the expansion \eqref{fipomb.26}, from which \eqref{fipomb.33} follows.
\end{proof}

The naturality of the definition of the cusp algebra means that for a
fibration with model fibre a compact manifold with boundary $Z$ we
may choose a smooth family of cusp structures, for instance by
choosing a cusp structure on the total space, and so define the bundle
of fibre-wise cusp operators $\Psi^m_{\cusp}(M/B;E,F)$ for any bundles
$E$ and $F$ over the total space and any $m.$ Then the short exact
sequence \eqref{CnvRcn.22} on each fibre becomes 
\begin{multline}
\Psi^{m-1}_{\cusp}(M/B;E,F)\longrightarrow\Psi^m_{\cusp}(M/B;E,F)\\
\overset{\sigma_m}\longrightarrow \CI({}^{\cusp}S^*(M/B);\hom(E,F)\otimes R^m).
\label{CnvRcn.24}\end{multline}
Since $R^m$ is trivial (canonically so if $m=0$) and
${}^{\cusp}S^*(M/B)$ is isomorphic to the `true' fibre cosphere bundle
$S^*(M/B),$ the image space in \eqref{CnvRcn.24} is an appropriate
replacement for the corresponding space of symbols in the non-boundary
case of Atiyah and Singer. In particular we may consider a family
$P\in \Psi^m_{\cusp}(M/B;E,F)$ to be `elliptic' if $\sigma_m(P)$ has
an inverse $s\in\CI({}^{\cusp}S^*(M/B);\hom(F,E)\otimes R^{-m});$
this is equivalent to pointwise invertibility of the
symbol. Ellipticity, as in the boundaryless case, is equivalent to the
existence of a parametrix $S\in\Psi^{-m}_{\cusp}(M/B;F,E)$ such that
$\Id-PS\in\Psi^{-\infty}_{\cusp}(M/B;E)$ and $\Id-SP\in
\Psi^{-\infty}_{\cusp}(M/B;F).$ However, the important difference, due
to the presence of the boundary, is that these errors are not
necessarily compact as families of operators on $L^2$ and
correspondingly, $P$ need not be a family of Fredholm operators.

\providecommand{\bysame}{\leavevmode\hbox to3em{\hrulefill}\thinspace}
\providecommand{\MR}{\relax\ifhmode\unskip\space\fi MR }
\providecommand{\MRhref}[2]{%
  \href{http://www.ams.org/mathscinet-getitem?mr=#1}{#2}
}
\providecommand{\href}[2]{#2}

\end{document}